\documentclass{amsart}
\usepackage[cp850]{inputenc}
\usepackage{graphics}
\usepackage{graphicx}
\usepackage{amssymb}

\newtheorem{theorem}{Theorem}
\newtheorem{lemma}[theorem]{Lemma}
\newtheorem{proposition}[theorem]{Proposition}

\theoremstyle{definition}

\newtheorem{example}[theorem]{Example}

\theoremstyle{remark}

\usepackage{amssymb,amsmath}

\def\eop{{\vrule height 6pt width 5pt depth 0pt}}\smallskip
\def\<{{\langle }}
\def\>{{\rangle }}

\newcommand{\g}[2]{\mbox{$\langle #1 ,#2 \rangle$}}
\newcommand{\s}[1]{\mbox{$\mathbb{S}^{#1}$}}
\newcommand{\R}[1]{\mbox{$\mathbb{R}^{#1}$}}

\newcommand{\rf}[1]{\mbox{(\ref{#1})}}
\newcommand{\rl}[1]{{~\ref{#1}}}

\def\indM{\mbox{$\mathrm{Ind}(M)$}}
\def\indT{\mbox{$\mathrm{Ind}_T(M)$}}

\def\beq{\begin{equation}}
\def\eeq{\end{equation}}

\begin{document}

\title[A characterization of quadric CMC hypersurfaces of spheres]
{A characterization of quadric constant mean curvature hypersurfaces of spheres}

\author{Luis J. Al\'\i as}
\address{Departamento de Matem\'{a}ticas, Universidad de Murcia, Campus de Espinardo,
E-30100 Espinardo, Murcia, Spain} \email{ljalias@um.es}
\thanks{L.J. Al\'\i as was partially supported by MEC project MTM2007-64504, and Fundaci\'{o}n S\'{e}neca
project 04540/GERM/06, Spain. This research is a result of the activity developed within the framework of the
Programme in Support of Excellence Groups of the Regi\'{o}n de Murcia, Spain, by
Fundaci\'{o}n S\'{e}neca, Regional Agency for Science and Technology (Regional Plan for
Science and Technology 2007-2010).}

\author{Aldir Brasil Jr.}
\address{Departamento de Matem\'{a}tica, Universidade Federal do Cear{\'{a}},  Campus do Pici,
60455-760 Fortaleza-Ce, Brazil} \email{aldir@mat.ufc.br}
\thanks{A. Brasil Jr. was partially supported by CNPq, Brazil, 306626/2007-1.}

\author{Oscar Perdomo }
\address{Department of Mathematical Sciences, Central Connecticut State University, New Britain, CT 06050, USA}
\email{perdomoosm@ccsu.edu}




\dedicatory{Dedicated to the memory of Professor Luis J. Al\'\i as-P\'{e}rez}


\begin{abstract}
Let $\phi:M\to\mathbb{S}^{n+1}\subset\mathbb{R}^{n+2}$  be an immersion of a complete $n$-dimensional oriented manifold. For any $v\in\mathbb{R}^{n+2}$, let us denote by $\ell_v:M\to\mathbb{R}$ the function given by
$\ell_v(x)=\<\phi(x),v\>$ and by $f_v:M\to\mathbb{R}$, the function given by $f_v(x)=\<\nu(x),v\>$, where $\nu:M\to\mathbb{S}^{n}$ is a Gauss map. We will prove that if $M$ has constant mean curvature, and, for some $v\ne{\bf 0}$ and
some real number $\lambda$, we have that $\ell_v=\lambda f_v$, then, $\phi(M)$ is either a totally umbilical sphere or a Clifford hypersurface. As an application, we will use this result to prove that the weak stability index of any
compact constant mean curvature hypersurface $M^n$ in $\mathbb{S}^{n+1}$ which is neither totally umbilical nor a Clifford hypersurface and has constant scalar curvature is greater than or equal to $2n+4$.
\end{abstract}

\maketitle

\section{Introduction}
Let $\phi:M\to\s{n+1}\subset\R{n+2}$ be an immersion of a complete $n$-dimensional
oriented manifold. For every $x\in M$ we will denote by $T_xM$ the tangent space of
$M$ at $x$. Sometimes, specially when we are dealing with local aspects of
$M$, we will identify $M$ with the set $\phi (M)\subset\R{n+2}$, and
the space $T_xM$ with the linear subspace $d\phi_x(T_xM)$ of $\R{n+2}$.
Let us denote by $\nu:M\to\s{n+1}\subset\R{n+2}$, a normal unit vector field along
$M$, i.e., for every $x\in M$, $\nu(x)$ is perpendicular to the vector $x$
and to the vector space $T_xM$. The shape operator $A_x:T_xM\to T_xM$, is
given by $A_x(v)=-d\nu_x(v)=-\beta^{\prime}(0)$ where
$\beta(t)=\nu(\alpha(t))$ and $\alpha(t)$ is any smooth curve in $M$ such
that $\alpha(0)=x$ and $\alpha^{\prime}(0)=v$. It can be shown that the
linear map $A_x:T_xM\to T_xM$ is symmetric, therefore it has $n$ real
eigenvalues $\kappa_1(x),\ldots,\kappa_{n}(x)$. These eigenvalues are
known as the principal curvatures of $M$ at $x$. The mean curvature of $M$
at $x$ is the average of the principal curvatures,
\[
H(x)=\frac{\kappa_1(x)+\cdots+\kappa_{n}(x)}{n},
\]
and the norm square of the shape operator is defined by the equation
\[
\|A\|^2(x)=\mathrm{trace}(A_x^2)=\kappa_1^2(x)+\cdots+\kappa_{n}^2(x).
\]

\subsection{Examples: Totally umbilical spheres and Clifford hypersurfaces} In
this section we will describe two families of examples that are related
with the main result of this paper.

\begin{example}
\label{spheres}
Let $v\in \R{n+2}$ be a  fixed unit
vector and $c$ a real number with $|c|<1$. Let us define
\[
\s{n}(v,c)=\{x\in \s{n+1}:\<x,v\>=c\}.
\]
Clearly, $\s{n}(v,c)$ is a hypersurface of $\s{n+1}$. In this case the map
$\nu:\s{n}(v,c)\to\s{n+1}$ given by
\[
\nu(x)=\frac{1}{\sqrt{1-c^2}}\left(v-cx\right)
\]
is a normal unit vector field along $\s{n}(v,c)$. Therefore, for every
$x\in\s{n}(v,c)$ the shape operator $A_x$ is the map $c(1-c^2)^{-{1\over 2}}I$, where $I$ is the
identity map, and
\[
\kappa_1(x)=\cdots=\kappa_{n}(x)=\frac{c}{\sqrt{1-c^2}}
\]
for all $x\in\s{n}(v,c)$. It is not difficult to show that these examples are
the only totally umbilical complete hypersurfaces of $\s{n+1}$. In this case
\[
H=\frac{c}{\sqrt{1-c^2}} \quad\hbox{and} \quad \|A\|^2=\frac{nc^2}{1-c^2}
\]
are both constant on $\s{n}(v,c)$.
\end{example}

\begin{example}
Given any integer $k\in \{1,\ldots,n-1\}$ and any real number
$r\in(0,1)$, let us define $\ell=n-k$ and
\begin{eqnarray*}
M_{k}(r) & = & \{(x,y)\in\R{k+1}\times\R{\ell+1}: \|x\|^2=r^2 \textrm{ and } \|y\|^2=1-r^2\}\\
{} & = & \s{k}(r)\times\s{n-k}(\sqrt{1-r^2})\subset\s{n+1}.
\end{eqnarray*}
It is not difficult to see that for any $(x,y)\in M_{k}(r)$ one gets
\[
T_{(x,y)}M_{k}(r) =\{(v,w)\in \R{k+1}\times\R{{\ell+1}}: \<x,v\>=0\quad \hbox{and}\quad \<w,y\>=0\}
\]
Therefore, the map $\nu:M_{k}(r)\to\s{n+1}$ given by
\[
\nu(x,y)=({\sqrt{1-r^2}\over r}x,-{r\over \sqrt{1-r^2}}y)
\]
defines a normal unit vector field along $M_{k}(r)$, i.e. it is a Gauss
map on $M_{k}(r)$. Notice that the vectors in $T_{(x,y)}M_{k}(r)$ of the
form $(v,{\bf 0})$ define a $k$ dimensional space. A direct
computation, using the expression for $\nu$, gives us that if
$(v,{\bf 0})\in T_{(x,y)}M_{k}(r)$, then,
\[
A_{(x,y)}(v,{\bf 0})=-\frac{\sqrt{1-r^2}}{r}(v,{\bf 0}).
\]
Therefore $-\sqrt{1-r^2}/r$ is an eigenvalue of $A_{(x,y)}$ with
multiplicity $k$. In the same way we can show that $r/\sqrt{1-r^2}$ is an eigenvalue of $A_{(x,y)}$
with multiplicity $\ell$. Therefore, the principal curvatures of $M_{k}(r)$ are given by
\[
\kappa_1(x,y)=\cdots=\kappa_{k}(x,y)=-\frac{\sqrt{1-r^2}}{r}, \quad
\kappa_{k+1}(x,y)=\cdots=\kappa_n(x,y)=\frac{r}{\sqrt{1-r^2}},
\]
and we also have that
\[
H=\frac{nr^2-k}{nr\sqrt{1-r^2}} \quad\hbox{and}\quad
\|A\|^2=\frac{k}{r^2}+\frac{n-k}{1-r^2}-n
\]
are both constant. Hypersurfaces that, up to a rigid motion, are equal to $M_{k}(r)$ for some $k$ and $r$, are called
Clifford hypersurfaces.
\end{example}

\subsection{Two families of geometric functions on hypersurfaces in spheres}
\label{ss1.2}
Given a fixed vector $v\in \R{n+2}$, let us define the functions
$\ell_v:M\to \R{}$ and $f_v:M\to \R{}$ by $\ell_v(x)=\<\phi(x),v\>$ and
$f_v(x)=\<\nu(x),v\>$, where $\nu:M\to\s{n+1}$ is a Gauss map. When we
consider all possible $v\in \R{n+2}$ we obtain the families
\[
V_1=\{\ell_v : v\in\R{n+2}\} \quad\hbox{and}\quad V_2=\{f_v : v\in\R{n+2}\}.
\]

These two families are very useful in the study of the spectrum of
important elliptic operators defined on $M$ like the Laplacian an the
stability operator. For example, in \cite{S1} and \cite{S2}, Solomon
computed the whole spectrum for the Laplace operator of every minimal
isoparametric hypersurface of degree 3 in spheres using these two families
of functions. For the totally umbilical spheres $\s{{n}}(v,c)$ we have
that if $c=0$, then dim$(V_1)=n+1$ and dim$(V_2)=1$. Indeed, it is not
difficult to prove that if for some compact hypersurface $M^n$ in \s{n+1},
we have that either $\hbox{dim}(V_1)<n+2$ or $\dim(V_2)<n+2$, then
$M=\s{{n}}(v,0)$ for some unit vector $v\in\R{n+2}$, \cite[Lemma 3.1]{P}.

If we take $c\ne 0$, and we consider the example $\s{n}(v,c)$ we observe
that if $w\in \R{n+2}$ is a vector perpendicular to the vector $v$, then
\[
f_w=-\frac{c}{\sqrt{1-c^2}}\ell_w.
\]
We also have this kind of relation between the function $f_w$ and the function $\ell_w$ in the Clifford
hypersurfaces; more precisely, if we consider the example $M_{k}(r)$ and
we take $w=(w_1,\dots,w_{k+1},0,\dots,0)\in \R{n+2}$ then we have that
\[
f_w=\frac{\sqrt{1-r^2}}{r}\ell_w.
\]
Also, if we take $w=(0,\dots,0,w_{k+2},\dots,w_{n+2})\in\R{n+2}$, then, we have that
\[
f_w=-\frac{r}{\sqrt{1-r^2}}\ell_w.
\]
In this paper we will prove that these two examples are the only hypersurfaces
with constant mean curvature in \s{n+1} where the relation $f_w=\lambda \ell_w$, for
some non-zero vector $w\in \R{n+2}$, is possible. More precisely, we
will prove the following result.
\begin{theorem}
\label{maintheorem}
Let $\phi:M\to\s{n+1}\subset\R{n+2}$ be an immersion with constant mean curvature of a complete $n$-dimensional oriented
manifold. If for some non-zero vector $v\ne{\bf 0}$ and some real number $\lambda$, we have that $\ell_v=\lambda f_v$,
then, $\phi(M)$ is either a totally umbilical sphere or a Clifford hypersurface.
\end{theorem}

Recall that
constant mean curvature hypersurfaces in \s{n+1} are characterized as critical points of the area functional
restricted to variations that preserve a certain volume function. As is well-known, the Jacobi operator of this
variational problem is given by $J=\Delta+\|A\|^2+n$, with associated quadratic form given by
\[
Q(f)=-\int_MfJf
\]
and acting on the space
\[
\mathcal{C}_T^\infty(M)=\{ f\in\mathcal{C}^\infty(M) : \mbox{$\int_M f=0$} \}.
\]
Precisely, the restriction $\int_M f=0$ means that the variation associated to $f$ is volume preserving.

In contrast to the case of minimal hypersurfaces, in the case of
hypersurfaces with constant mean curvature
one can consider two different eigenvalue problems: the usual Dirichlet problem, associated with the quadratic form $Q$
acting on the whole space of smooth functions on $M^n$, and the so called \textit{twisted} Dirichlet problem, associated
with the same quadratic form $Q$, but restricted to the subspace of smooth functions satisfying the additional condition $\int_Mf=0$.
Similarly, there are two different notions of stability and
index, the \textit{strong stability} and \textit{strong index}, denoted by \indM\ and associated
to the usual Dirichlet problem, and the \textit{weak stability} and \textit{weak index}, denoted
by \indT\ and associated to the twisted Dirichlet problem.
Specifically, the strong index of the hypersurface is characterized as
\[
\mathrm{Ind}(M)=
\max\{ \mathrm{dim}V : V\leqslant\mathcal{C}^\infty(M), \quad Q(f)<0 \quad \mbox{for every } f\in V \},
\]
and $M$ is called strongly stable if and only if $\mathrm{Ind}(M)=0$. On the other hand, the weak stability index of $M^n$ is characterized by
\[
\mathrm{Ind}_T(M)=
\max\{ \mathrm{dim}V : V\leqslant\mathcal{C}_T^\infty(M), \quad Q(f)<0 \quad \mbox{for every } f\in V \},
\]
and $M$ is called weakly stable if and only if $\mathrm{Ind}_T(M)=0$. From a geometrical point of view, the weak index is
more natural than the strong index. However, from an analytical point of view, the strong index is more natural and
easier to use (for further details, see \cite{A}).

As an application of our Theorem \ref{maintheorem}, we will prove that the weak stability index of a compact
constant mean curvature hypersurface $M^n$ in $\s{n+1}$ with constant scalar curvature must be greater than
or equal to $2n+4$ whenever $M^n$ is neither a totally umbilical sphere nor a Clifford hypersurface (see Theorem \ref{theorindex}).
This result complements the one obtained in \cite{ABP} where the authors showed that
the weak index of a compact constant mean curvature hypersurface $M^n$ in \s{n+1} which
is not totally umbilical and has constant scalar curvature is greater than or equal to $n+2$,
with equality if and only if $M^n$ is a Clifford hypersurface
$M_k(r)=\mathbb{S}^{k}(r)\times\mathbb{S}^{n-k}(\sqrt{1-r^2})$ with radius
$\sqrt{k/(n+2)}\leqslant r\leqslant\sqrt{(k+2)/(n+2)}$. At this respect, it is
worth pointing out that the weak stability index of the Clifford hypersurfaces $M_k(r)$ depends on $r$, reaching its
minimum value $n+2$ when $\sqrt{k/(n+2)}\leqslant r\leqslant\sqrt{(k+2)/(n+2)}$, and converging to $+\infty$ as $r$ converges either
to $0$ or $1$ (see \cite[Section 3]{ABP} for further details).

\section{Preliminaries and auxiliary results}
Let us start this section by  computing the gradient of the functions $\ell_v$
and $f_v$.  For any fixed vector in $\R{n+2}$, let us define the tangent
vector field $v^{\top}:M\to\R{n+2}$ by
\[
v^{\top}(x)=v-\ell_v(x)x-f_v(x)\nu(x)\qquad \hbox{for all $x\in M$},
\]
where, as in the previous section, $\nu:M\to\R{n+2}$ is a Gauss map. Clearly, $v^{\top}$ is a tangent vector field
on $M$ because $\<v^{\top}(x),x\>=0$ and $\<v^{\top}(x),\nu(x)\>=0$ for every $x\in M$. More precisely, $v^{\top}(x)$ is
the orthogonal projection of the vector $v$ on $T_xM$.
\begin{proposition}
If $M^n$ is a smooth hypersurface of $\s{n+1}$
and $A$ denotes its shape operator with respect to the unit normal vector
field $\nu:M\to \R{n+2}$ then, the gradient of the functions $\ell_v$ and
$f_v$ are given by:
\[
\nabla \ell_v = v^{\top}, \quad  \quad \nabla f_v = -A(v^{\top}).
\]
\end{proposition}
\begin{proof}
For any vector $w\in T_xM$, let $\alpha:(-\varepsilon,\varepsilon)\to M$ be a curve such that
$\alpha(0)=x$ and $\alpha^{\prime}(0)=w$. Notice that
\[
d\ell_v(w)={d\ell_v(\alpha(t))\over dt}\big{|}_{t=0}={d\<\alpha(t),v\>\over dt}\big{|}_{t=0}=
\<\alpha^{\prime}(0),v\>=\<w,v^{\top}(x)\>.
\]
Since the equality above holds true for every $w\in T_xM$ and $v^{\top}(x)\in T_xM$, then, $\nabla \ell_v(x)=v^{\top}(x)$.
For the function $f_v$, we have
\begin{eqnarray*}
df_v(w) & = & {df_v(\alpha(t))\over dt}\big{|}_{t=0}={\<\nu(\alpha(t)),v\>\over dt}\big{|}_{t=0}=
\<d\nu(\alpha^{\prime}(0)),v\>\\
{} & = & -\<A(w),v^{\top}(x)\>=\<w,-A(v^{\top}(x))\>.
\end{eqnarray*}
Therefore, $\nabla f_v(x)=-A(v^{\top}(x))$.
\end{proof}

We also have the following expressions for the Laplacian of the functions $\ell_v$ and $f_v$.
\begin{proposition}
\label{prop5}
If $M^n$ is a smooth hypersurface of $\s{n+1}$ with constant mean curvature $H$, and $A$ denotes the shape operator with
respect to the unit normal vector field $\nu:M\to \R{n+2}$ then, the Laplacian of the functions $\ell_v$ and $f_v$ are
given by:
\[
\Delta \ell_v = -n\ell_v+nH f_v, \quad  \quad
\Delta f_v = -\|A\|^2f_v+nH\ell_v.
\]
\end{proposition}
\begin{proof}
For any vector $w\in T_xM$, we have
\[
\nabla_w\nabla\ell_v=\nabla_wv^\top=-\ell_v(x)w+f_v(x)A_x(w),
\]
where $\nabla$ denotes here the intrinsic derivative on $M$. Let $\{ e_1,\ldots, e_n\}$ be an orthonormal basis of
$T_xM$. Then, the Laplacian of $\ell_v$ at the point $x$ is given by
\[
\Delta\ell_v(x)=\sum_{i=1}^n\g{\nabla_{e_i}\nabla\ell_v}{e_i}=-n\ell_v(x)+\mathrm{tr}(A_x)f_v(x)=
-n\ell_v(x)+nHf_v(x).
\]
On the other hand, using Codazzi equation we also have that
\begin{eqnarray*}
\nabla_w\nabla f_v & = & -\nabla_w(A(v^\top))=-(\nabla_wA)(v^\top(x))-A_x(\nabla_wv^\top)\\
{} & = & -(\nabla_{v^\top(x)}A)(w)+\ell_v(x)A_x(w)-f_v(x)A_x^2(w).
\end{eqnarray*}
Therefore
\begin{eqnarray*}
\Delta f_v (x)& = & \sum_{i=1}^n\g{\nabla_{e_i}\nabla f_v}{e_i} \\
{} & = & -\sum_{i=1}^n\g{(\nabla_{v^\top(x)}A)(e_i)}{e_i}+nH\ell_v(x)-\|A\|^2(x)f_v(x)\\
{} & = & -n\g{v^\top(x)}{\nabla H(x)}+nH\ell_v(x)-\|A\|^2(x)f_v(x)\\
{} & = & nH\ell_v(x)-\|A\|^2(x)f_v(x),
\end{eqnarray*}
since the mean curvature $H$ is constant.
\end{proof}

The following two lemmas will be used in the proof of our main theorem. The first one is an elementary geometric lemma whose proof is left to the reader.
\begin{lemma}
\label{lemma2.1}
Let $M^n$ be a smooth hypersurface of \s{n+1} and let $\alpha:I\subset\R{}\to M$ be a regular curve such that
\[
\alpha^{\prime\prime}(t)=f(t)\alpha^{\prime}(t)+\eta(t)
\]
where $f:I\to \R{}$ is a smooth function and $\eta:I\to \R{n+2}$
is a normal vector field along $\alpha$, i.e. $\eta(t)$ is orthogonal to $T_{\alpha(t)}M$.
If $s=s(t)$ is the arc-length parameter for $\alpha$, then
$\beta(s)=\alpha(t(s))$ satisfies that $\beta^{\prime\prime}(s)$ is a
normal vector field along $\beta$, i.e. $\beta$ is a geodesic in $M$.
\end{lemma}

The other one is an algebraic lemma.
\begin{lemma}
\label{lemma2.2}
If $p_1(X)=b_1X+c_1,\ldots, p_k(X)=b_kX+c_k$ are $k$ polynomials of degree 1, $k\geq 2$, with the property that
$c_i/b_i\ne c_j/b_j$ whenever $i\ne j$, then, the polynomials
\[
q_i=\Pi_{j=1,j\ne i}^k p_j
\]
are linearly independent. Moreover, an equation of the form
\[
{a_1\over p_1(X)}+\cdots +{a_k\over p_k(X)}= d
\]
with $a_i$ and $d$ real numbers, can not hold true unless all the $a_i$'s and $d$ are zero.
\end{lemma}
\begin{proof}
By the condition on the numbers $c_j/b_j$ we
have that at $X_i=-c_i/b_i$ every polynomial $q_j$, except
the polynomial $q_i$, vanishes. Therefore, if there exists constants $\alpha_i$ such that
\[
\alpha_1q_1(X)+\cdots+ \alpha_kq_k(X)=0
\]
then, taking $X=X_i$ we get that $\alpha_i=0$ for every $i$.  Therefore, the
polynomials $q_i$'s are linearly independent. On the other hand,
notice that the second equation in the lemma can be written as
\[
a_1q_1(X)+\cdots+ a_kq_k(X)=dR(X)
\]
where $R$ is a polynomial of degree $k$. Since the expression on
the left of the last equation is a polynomial of degree $k-1$, we
obtain that the constant on the right hand side must be zero. Then
the second part of the lemma follows by the independence of the
polynomials $q_i$'s.
\end{proof}

\section{Proof of Theorem\rl{maintheorem}}
We are now ready to give our main argument and prove Theorem\rl{maintheorem}.
Since most of the arguments are local and the thesis of the
theorem is on $\phi(M)$ and not on $M$, we will identify $M$ with
$\phi(M)$ and $T_xM$ with $T_{\phi(x)}M$. By multiplying the equation
$\ell_v=\lambda f_v$ by an appropriated constant we may assume that $|v|=1$.
We will also assume that $\ell_v$ is not constant, otherwise
$\phi(M)\subset\s{n}(v,c)$ for some $c$, which implies, using the completeness of $M$,
that $\phi(M)=\s{n}(v,c)$.

Notice that, since $\ell_v$ is not constant, then $\lambda\ne 0$.
Taking the gradient in both sides of the expression $\ell_v=\lambda f_v$ we obtain that
\beq
\label{2.1}
A(v^{\top}(x))= -\lambda^{-1}v^{\top}(x)
\eeq
at every point $x\in M$.

\

\noindent{\bf Step 1: The integral curves of $v^\top$ in $M$ are Euclidean circles.}
Let us take a point $x\in M$ such that $\nabla\ell_v(x)=v^{\top}(x)$ does not vanish.
Let $\alpha_x(t)$ be the integral curve of the vector field $v^{\top}$ such that $\alpha_x(0)=x$. Since
\begin{eqnarray*}
\alpha_x^{\prime}(t) & = & v^{\top}(\alpha_x(t))=v-\ell_v(\alpha_x(t))\alpha_x(t)-f_v(\alpha_x(t))\nu(\alpha_x(t))\\
{} & = &  v-\ell_v(\alpha_x(t))\big{(}\alpha_x(t)+\lambda^{-1}\nu(\alpha_x(t))\big{)}
\end{eqnarray*}
then,
\begin{eqnarray}
\label{2.2}
\nonumber \alpha_x^{\prime\prime}(t) & = &
-\<\nabla\ell_v(\alpha_x(t)),\alpha_x^{\prime}(t)\>\big{(}\alpha_x(t)+\lambda^{-1}\nu(\alpha_x(t))\big{)}\\
\nonumber {} & {} & -\ell_v(\alpha_x(t))\big{(}\alpha_x^{\prime}(t)-\lambda^{-1}A(\alpha_x^{\prime}(t))\big{)}\\
{} & = & -|\alpha_x^{\prime}(t)|^2 \big{(}\alpha_x(t)+\lambda^{-1}\nu(\alpha_x(t))\big{)}
-\ell_v(\alpha_x(t))\big{(}\alpha_x^{\prime}(t) -\lambda^{-1}A(v^{\top}(\alpha_x(t)))\big{)}\\
\nonumber {} & = & -|\alpha_x^{\prime}(t)|^2\big{(}\alpha_x(t)+\lambda^{-1}\nu(\alpha_x(t))\big{)}
-\ell_v(\alpha_x(t))\big{(}\alpha_x^{\prime}(t) +\lambda^{-2}v^{\top}(\alpha_x(t))\big{)}\\
\nonumber {} & = & f(t)\alpha_x^{\prime}(t)+\eta(t).
\end{eqnarray}
Here
\beq
\label{1.2}
\eta(t)=-|\alpha_x^{\prime}(t)|^2\big{(}\alpha_x(t)+\lambda^{-1}\nu(\alpha_x(t))\big{)}
\eeq
is a normal vector field along $\alpha_x$ and
\[
f(t)=-\big{(}1+\lambda^{-2}\big{)}\ell_v(\alpha_x(t)).
\]
Therefore if $s=s(t)$ is the arc-length
parameter for the curve $\alpha_x$ with $s(0)=0$, and $t=t(s)$ is
the inverse of the function $s=s(t)$, we have, by Lemma\rl{lemma2.1},  that
$\beta_x(s)=\alpha_x(t(s))$ is a geodesic in $M$. Moreover, from \rf{2.2} and \rf{1.2} we also get that
\beq
\label{2.3}
\beta_x^{\prime\prime}(s)=\frac{1}{|\alpha_x^{\prime}(t(s))|^2}\eta(t(s))=-\beta_x(s)-\lambda^{-1}\nu(\beta_x(s)).
\eeq

If we differentiate \rf{2.3}, we get that the
function $\beta_x^{\prime}$ moves along a circle because it
satisfies the equation
\[
(\beta_x^{\prime})^{\prime\prime}(s)+(1+\lambda^{-2})\beta_x^{\prime}(s)=0.
\]
More precisely, if we define $w=\sqrt{1+\lambda^{-2}}>0$, then
\[
\beta_x^{\prime}(s)=\beta_x^{\prime}(0)\cos{(ws)}+w^{-1}\beta_x^{\prime\prime}(0)\sin{(ws)},
\]
and
\beq
\label{2.6}
\beta_x(s)=w^{-1}\beta_x^{\prime}(0)\sin{(ws)}-w^{-2}\beta_x^{\prime\prime}(0)\cos{(ws)}+
\beta_x(0)+ w^{-2}\beta_x^{\prime\prime}(0).
\eeq

\

\noindent{\bf Step 2: The intersection $N=M\cap\mathbb{S}^n(v,0)$ is non-empty.}
Let us compute $\beta_x^{\prime}(0)$ and $\beta_x^{\prime\prime}(0)$
in order to obtain an explicit expression for $\beta_x^{\prime}(s)$.
From the definition of $\beta_x$ we have that $\beta_x(0)=x$ and

\beq
\label{2.7}
\beta_x^{\prime}(0)=\frac{\alpha_x^\prime(0)}{|\alpha_x^\prime(0)|}=
\frac{v^{\top}(x)}{|v^{\top}(x)|}.
\eeq
Notice that
\[
v^{\top}(y)=v-\ell_v(y)y-f_v(y)\nu(y)=v-\ell_v(y)y-\lambda^{-1}\ell_v(y)\nu(y)
\]
at every point $y\in M$. Therefore
\[
|v^{\top}(y)|^2= 1-\ell_v(y)^2-\lambda^{-2}\ell_v(y)^2=1-w^2\ell_v(y)^2. 
\]
From this last expression we obtain that $-w^{-1}\le \ell_v(y)\le w^{-1}$, at every $y\in M$, and
\beq
\label{2.9}
v^{\top}(y)={\bf 0} \mbox{ if and only if } \ell_v(y)=\pm w^{-1}.
\eeq

Let us define $a=\ell_v(x)$, and $b=\sqrt{w^{-2}-a^2}$. By \rf{2.9} we
have that $b>0$, because $\nabla \ell_v(x)=v^{\top}(x)\ne {\bf 0}$.
With this notation, we obtain that $|v^{\top}(x)|^2=1-w^2a^2=w^2b^2$, and
\begin{eqnarray*}
\<\beta_x(0),v\>& = & \ell_v(x)=a\\
\<\beta_x^\prime (0),v\> & = & \<{v^{\top}(x)\over |v^{\top}(x)|},v\>
=\<{v^{\top}(x)\over |v^{\top}(x)|},v^{\top}(x)\>
=|v^{\top}(x)|=\sqrt{1-w^2a^2}=wb \\
\<\beta_x^{\prime\prime}(0),v\> & = & \<-\beta_x(0)-\lambda^{-1}\nu(\beta_x(0)),v\>=
-a-\lambda^{-2}a=-w^2a,
\end{eqnarray*}
where we have used \rf{2.3} to derive the last equation. Now, using these equations jointly with \rf{2.6} we
get that
\[
\ell_v(\beta_x(s))=\<\beta_x(s),v\>=a\cos{(ws)}+b\sin{(ws)}.
\]
Notice that $(wa)^2+(wb)^2=1$ with $wb>0$. Therefore for some $s_1\in (-{\pi\over 2w},{\pi\over 2w})$ we have
\[
-wa=\sin{(ws_1)}\quad\hbox{and}\quad wb=\cos{(ws_1)},
\]
so that
\[
\ell_v(\beta_x(s))=a\cos{(ws)}+b\sin{(ws)}=w^{-1}\sin{(ws-ws_1)}.
\]
Notice that when $s$ moves from $0$ to $s_1$, we have that
$\ell_v(\beta_x(s))$ never reaches the values $\pm w^{-1}$, therefore by \rf{2.9}
$v^\top(\beta_x(s))\neq{\bf 0}$ and all these $\beta_x(s)$ belong to the integral curve of the
vector field $v^{\top}$.
In particular, $\ell_v(\beta_x(s_1))=0$ and $v^{\top}(\beta_x(s_1))=v\neq{\bf 0}$.
This argument shows that
\[
N=\ell_v^{-1}(0)=\{ y\in M : \ell_v(y)=0 \}
\]
is not empty. Observe that if we were assuming that $M$ were compact instead of
complete, the fact that $N=\ell_v^{-1}(0)$ is not empty would have
followed from the fact that the function $\ell_v$ must reach its
maximum value and a minimum value on $M$, and the fact that necessarily
these values must be $\pm w^{-1}$, since $\nabla\ell_v=v^\top$ must vanish at its critical points.
From now on we will assume that the $x$ that we were considering before is an element in $N$, i.e,
we will assume that $a=0$, and therefore $b=w^{-1}$ and $s_1=0$.

\

\noindent{\bf Step 3: The intersection $N=M\cap\mathbb{S}^n(v,0)$ as a hypersurface of $M$ and as a hypersurface of $\mathbb{S}^n(v,0)$.}
Clearly the set $N\subset M^n$ is an $(n-1)$-dimensional manifold because $0$
is a regular value of the function $\ell_{v}$ on $M$. Moreover, for every
$x\in N$ we have that $\nabla \ell_v(x)=v^\top(x)=v$ is a constant vector, and therefore $N$ is a totally geodesic
hypersurface of $M$.
Notice that for every $x\in N$ we have that $v\in T_xM$ and
$A_x(v)=-\lambda^{-1}v$. Therefore we can take vectors
$v_1,\dots,v_{n-1}$ in $T_xM$, all of them orthogonal to $v$, such
that $A_x(v_i)=\lambda_i(x)v_i$. Since the vectors $v_i$'s are
perpendicular to $v=\nabla\ell_v(x)$, they form a basis for $T_xN$.
On the other hand, notice that $N$ is also a hypersurface of the unit $n$-dimensional sphere
$\s{n}(v,0)$, and that for every $x\in N$, $\nu(x)$ gives a unit vector field normal to $N$ in $\s{n}(v,0)$ (see Figure 1).

\

\begin{center}
\includegraphics[width=10cm]{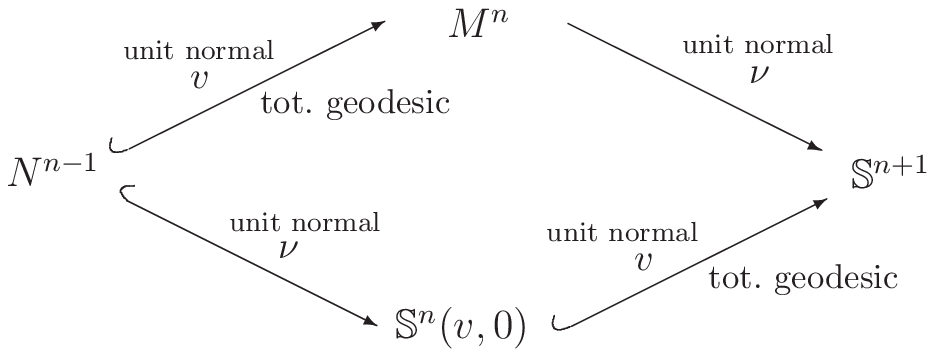}
\end{center}

\

Taking into account that $N$ is totally geodesic in $M^n$ and that $\s{n}(v,0)$ is totally geodesic in \s{n+1}, it
follows from the fact that $\nu$ is both normal to $M^n$ in \s{n+1} and normal to $N$ in $\s{n}(v,0)$ that, for every
$x\in N$, $\lambda_1(x),\ldots,\lambda_{n-1}(x)$ are the principal curvatures of $N$ as a hypersurface
of $\s{n}(v,0)$ with respect to $\nu$

\

\noindent{\bf Step 4: Computation of the principal curvatures of $M$ along the integral curves of $v^\top$.}
Under the assumption that $x\in N$, we obtain from \rf{2.7} that
\[
\beta_x^{\prime}(0)=v.
\]
Therefore, from \rf{2.3} and \rf{2.6} we get the following  expression for $\beta_x(s)$,
\beq
\label{2.12}
\beta_x(s)=w^{-1}\sin{(ws)}v+w^{-2}(\cos{(ws)}-1)(x+\lambda^{-1}\nu(x))+x.
\eeq
By differentiating two times this equation, and using the
equation \rf{2.3}, we obtain the following expression,
\beq
\label{2.13}
\nu(\beta_x(s))= \lambda w\sin(ws)v+\lambda\cos(ws)(x+\lambda^{-1}\nu(x))
-\lambda \beta_x(s).
\eeq
Recall that, if $s\in (-{\pi\over 2w},{\pi\over 2w})$, then
\beq
\label{2.13b}
|\ell_v(\beta_x(s))|<w^{-1} \quad \textrm{ and } \quad v^\top(\beta_x(s))\neq {\bf 0}.
\eeq

Observe that if $\gamma(t)$ is a smooth curve in $N$ such that $\gamma(0)=x$ and $\gamma^{\prime}(0)=v_i$, then by
\rf{2.12}, we have that the curve
\[
\gamma_s(t)=\beta_{\gamma(t)}(s)=
w^{-1}\sin{(ws)}v+w^{-2}(\cos{(ws)}-1)(\gamma(t)+\lambda^{-1}\nu(\gamma(t)))+\gamma(t)
\]
is a curve on $M$ such that $\gamma_s(0)=\beta_x(s)$. A direct computation shows that
\beq
\label{2.14}
\gamma_s^{\prime}(0)=
w^{-2}(\cos{(ws)}-1)(v_i-\lambda^{-1}\lambda_i(x)v_i)+v_i=\mu_i(x)v_i,
\eeq
where
\[
\mu_i(x)=\frac{\lambda(\lambda-\lambda_i(x))\cos(ws)+(1+\lambda\lambda_i(x))}{1+\lambda^2}.
\]
The computation above shows us that the vectors $v_i$'s are also elements in $T_{\beta_x(s)}M$ for every
$s\in (-{\pi\over 2w},{\pi\over 2w})$. Actually, it follows directly from \rf{2.14} that if
$\mu_i(x)\neq 0$
then $v_i=\gamma_s^{\prime}(0)/\mu_i(x)\in T_{\beta_x(s)}M$; hence by a continuity argument, since the
equation $\mu_i(x)=0$
has finitely many solutions on
$(-{\pi\over 2w},{\pi\over 2w})$, we conclude that $v_i\in T_{\beta_x(s)}M$ for every
$s\in (-{\pi\over 2w},{\pi\over 2w})$.

Recall that, by \rf{2.13b} and \rf{2.1}, $-\lambda^{-1}$ is a principal curvature at the point
$\beta_x(s)$, for every $s\in (-{\pi\over 2w},{\pi\over 2w})$, with associated principal direction in the direction of
$v^\top(\beta_x(s))\neq\mathbf{0}$.
Let us compute now the other $n-1$ principal curvatures of $M$ at the point $\beta_x(s)$. Since $\gamma(t)\in N$
for every $t$, then the expression \rf{2.13} holds true when replacing $x$ by $\gamma(t)$ and then we have that
\[
\nu(\gamma_s(t))=
\nu(\beta_{\gamma(t)}(s))=
\lambda w\sin(ws)v+\lambda\cos(ws)(\gamma(t)+\lambda^{-1}\nu(\gamma(t)))
-\lambda\gamma_s(t).
\]
Differentiating this equation with respect to $t$ at $t=0$ and using \rf{2.14}, we get that
\[
A_{\beta_x(s)}(\gamma'_s(0))=\mu_i(x)A_{\beta_x(s)}(v_i)=
-d\nu(\gamma_s^\prime(0))=(\lambda_i(x)-\lambda)\cos(ws)v_i+\lambda\mu_i(x)v_i.
\]
That is,
\[
A_{\beta_x(s)}(v_i)=
\left(\lambda+\frac{(\lambda_i(x)-\lambda)(1+\lambda^2)\cos(ws)}{\lambda(\lambda-\lambda_i(x))\cos(ws)+(1+\lambda\lambda_i(x))}\right)v_i
\]
Therefore, we get the following expression for the other $n-1$ principal curvatures at $\beta_x(s)$,
\begin{eqnarray}
\label{2.15}
\nonumber \lambda_i(\beta_x(s)) & = & \lambda+
{(\lambda_i(x)-\lambda)(1+\lambda^2)\cos(ws)\over\lambda(\lambda-\lambda_i(x))\cos(ws)+(1+\lambda\lambda_i(x))}\\
{} & = &
-\lambda^{-1}+
{(1+\lambda^2)(\lambda^{-1}+\lambda_i(x))\over\lambda(\lambda-\lambda_i(x))\cos{(ws)} +(1+\lambda\lambda_i(x))}.
\end{eqnarray}
Notice that, as it is supposed to be, when $s=0$, i.e at the point
$x$, the expression \rf{2.15} above reduces to $\lambda_i(x)$. Also
notice that if $\lambda_i(x)=-\lambda^{-1}$ then, the expression
\rf{2.15} reduces to $-\lambda^{-1}$ for every $s$.

\

\noindent{\bf Step 5: $M$ is isoparametric with at most two distinct principal curvatures.}
Now, we will use the hypothesis on the mean curvature of $M$. By \rf{2.15}, for every point $x\in N$ and every
$s\in (-{\pi\over 2w},{\pi\over 2w})$ we have that
\begin{eqnarray*}
nH & = & nH(\beta_x(s))=-\lambda^{-1}+\sum_{i=1}^{n-1}\lambda_i(\beta_x(s))\\
{} & = & -n\lambda^{-1}+(1+\lambda^2)
\sum_{i=1}^{n-1}\frac{\lambda^{-1}+\lambda_i(x)}{\lambda(\lambda-\lambda_i(x))\cos{(ws)}+(1+\lambda\lambda_i(x))}.
\end{eqnarray*}
That is,
\beq
\label{nH1}
\sum_{i=1}^{n-1}\frac{\lambda^{-1}+\lambda_i(x)}{\lambda(\lambda-\lambda_i(x))\cos{(ws)}+(1+\lambda\lambda_i(x))}=
\frac{n(H+\lambda^{-1})}{1+\lambda^2}.
\eeq
For every $x\in N$, let
\begin{eqnarray*}
I_1(x) & = & \{ i\in\{1,\ldots,n-1\} : \lambda_i(x)=-\lambda^{-1} \},\\
I_2(x) & = & \{ i\in\{1,\ldots,n-1\} : \lambda_i(x)=\lambda \},\\
I_3(x) & = & \{1,\ldots,n-1\}\setminus(I_1(x)\cup I_2(x)).
\end{eqnarray*}
Then \rf{nH1} can be written as
\beq
\label{nH2}
\sum_{i\in I_3(x)}\frac{\lambda^{-1}+\lambda_i(x)}{\lambda(\lambda-\lambda_i(x))\cos{(ws)}+(1+\lambda\lambda_i(x))}=d(x)
\eeq
where
\[
d(x)=\frac{n(H+\lambda^{-1})-n_2(x)(\lambda+\lambda^{-1})}{1+\lambda^2}.
\]
and $n_i(x)=\mathrm{card}(I_i(x))$. We claim that $I_3(x)=\emptyset$. Otherwise, for every $i\in I_3(x)$ let
$a_i(x)=\lambda^{-1}+\lambda_i(x)\neq 0$, $b_i(x)=\lambda(\lambda-\lambda_i(x))\neq 0$, and
$c_i(x)=1+\lambda\lambda_i(x)\neq 0$.  Thus, equation \rf{nH2} means that, for every $s\in (-{\pi\over 2w},{\pi\over 2w})$,
$\cos{(ws)}$ is a root of the polynomial equation on $X$
\beq
\label{nH3}
\sum_{i\in I_3(x)}\frac{a_i(x)}{b_i(x)X+c_i(x)}=d(x).
\eeq
If $\lambda_i(x)=\lambda_j(x)$  for every $i,j\in I_3(x)$ (in particular, if $n_3(x)=1$), then \rf{nH3}
becomes
\[
\frac{n_3(x)a_i(x)}{b_i(x)X+c_i(x)}=d(x),
\]
which can hold only if $a_i(x)=d(x)=0$. But this is a contradiction because $a_i(x)\neq 0$. Therefore, we can decompose
\[
I_3(x)=\bigcup_{i=1}^kJ_i(x), \quad k\geq 2,
\]
with $\lambda_{j_1}(x)=\lambda_{j_2}(x)$ if and only if $j_1,j_2\in J_i(x)$ for some $i$. In that case,
let $\lambda_{i}(x)=\lambda_{j}(x)$ for every $j\in J_i(x)$, and \rf{nH3} becomes
\beq
\label{nH4}
\sum_{i=1}^{k}\frac{m_i(x)a_i(x)}{b_i(x)X+c_i(x)}=d(x)
\eeq
with $m_i(x)=\mathrm{card}(J_i(x))>0$, $m_i(x)a_i(x)\neq 0$. But this contradicts our Lemma\rl{lemma2.2}, because
\[
\frac{c_i(x)}{b_i(x)}=\frac{1+\lambda\lambda_i(x)}{\lambda(\lambda-\lambda_i(x))}\neq
\frac{1+\lambda\lambda_j(x)}{\lambda(\lambda-\lambda_j(x))}=\frac{c_j(x)}{b_j(x)}
\]
for every $i\neq j$, $1\leq i,j\leq k$.

Summing up, $I_3(x)=\emptyset$ for every $x\in N$, which means that all the principal curvatures of $M$ at the points
of $N$ are constant and they are equal to either $-\lambda^{-1}$ or $\lambda$. From the expression \rf{2.15}, the same
happens along the geodesics $\beta_x(s)$ for every $s\in (-{\pi\over 2w},{\pi\over 2w})$. Taking into account that every
point of $M$ which is not a critical point of $\ell_v$ can be reached through a geodesic $\beta_x(s)$, we conclude that
the principal curvatures of $M$ are constant on the whole $M$ and they are equal to either $-\lambda^{-1}$ or $\lambda$.
That is, $M$ is a complete isoparametric hypersurface of \s{n+1} with at most two distinct principal curvatures, and from
the well known rigidity result by Cartan \cite{Ca} (see also \cite[Chaper 3]{CR}) we conclude that $M$ is either a totally
umbilical sphere (in the case that all its principal curvatures are equal to  $-\lambda^{-1}$) or it is either Clifford hypersurface
of the form $M_k(r)=\mathbb{S}^{k}(r)\times\mathbb{S}^{n-k}(\sqrt{1-r^2})$ with radius $0<r<1$ (in the case that the
principal curvatures  take both values).

This finishes the proof of Theorem\rl{maintheorem}.

\

\

Let us exhibit an example that shows that the condition on the
mean curvature to be constant is necessary in the previous result.
\begin{example}
Let $e_1=(1,0,\dots,0)\in\R{n+1}$ and $c=4/5$. From Example\rl{spheres} we know that the principal curvatures of
$\s{n-1}(e_1,c)\subset\s{n}$ are all equal to $4/3$. By perturbing $\s{n-1}(e_1,c)$ we
can find a hypersurface $N\subset\s{n}$ whose mean curvature is not constant and such that
all its principal curvatures $\lambda_i$ satisfy that
\beq
\label{2.17}
1<\lambda_i(x)<2 \quad \hbox{for every $x\in N$ and $i=1,\dots,n-1$}
\eeq
Let $M^n=\s{1}\times N$ and $\phi:M\to\s{n+1}\subset\R{n+2}$ the map given by
\[
\phi((\cos s,\sin s),x)= ({1\over\sqrt{2}}\sin(\sqrt{2}s),
{1\over 2}(x+\nu(x))\cos(\sqrt{2}s)+{1\over 2}(x-\nu(x))),
\]
where $x\in N\subset\s{n}\subset\R{n+1}$ denotes the points in $N$ and
$\nu:N\to\s{n}\subset\R{n+1}$ is a Gauss map of $N$. In particular, $\g{x}{\nu(x)}=0$.

Let ${\partial \over\partial s}=(-\sin s,\cos s)$ and
let $v_1,\dots,v_{n-1}$ be a basis of $T_xN$ such that
$-d\nu_x(v_i)=\lambda_i(x) v_i$.
Notice that $\bar{{\partial \over \partial s}}=((-\sin s,\cos s),{\bf 0})\in \R{{n+3}}$ and
$\bar{v}_1=(0,0,v_1),\dots,\bar{v}_{n-2}=(0,0,v_{n-2})$ form a basis for the tangent space of $M$ at
$p=((\cos s, \sin s),x)$. A direct computation shows that
\[
d\phi_p(\bar{{\partial \over \partial s}})=
(\cos(\sqrt{2}s),-{1\over \sqrt{2}}(x+\nu(x))\sin(\sqrt{2}s))
\]
and
\[
d\phi_p(\bar{v}_i)= {1\over 2}(0,((1-\lambda_i(x))\cos(\sqrt{2}s)+1+\lambda_i(x))v_i).
\]
By \rf{2.17}, the expression $(1-\lambda_i(x))\cos(\sqrt{2}s)+(1+\lambda_i(x))$ never vanishes,
therefore $\phi$ is an immersion. Moreover, it is easy to check that
$\tilde{\nu}:M\to \s{n+1}\subset\R{n+2}$ given by
\[
\tilde{\nu}(p)=({1\over \sqrt{2}}\sin(\sqrt{2}s),{1\over 2}(x+\nu(x))\cos(\sqrt{2}s)-{1\over 2}(x-\nu(x)))
\]
is a Gauss map on $M$. Using the expression for $\phi$ and for
$\tilde{\nu}$ we get that $\ell_v=f_v$ for $v=(1,0,\dots,0)\in \R{n+2}$.
\end{example}

\section{Stability index of hypersurfaces with constant mean curvature}
\label{sindex}
In this section, and as an application of our Theorem\rl{maintheorem}, we will prove that the weak stability index of a
compact constant mean curvature hypersurface $M^n$ in $\s{n+1}$ with constant scalar curvature must be greater than
or equal to $2n+4$ whenever $M^n$ is neither a totally umbilical sphere nor a Clifford hypersurface. Recall that
constant mean curvature hypersurfaces in \s{n+1} are critical points of the area functional
restricted to variations that preserve a certain volume function. The Jacobi operator of this
variational problem is given by $J=\Delta+\|A\|^2+n$, with associated quadratic form given by
\[
Q(f)=-\int_MfJf
\]
and acting on the space
\[
\mathcal{C}_T^\infty(M)=\{ f\in\mathcal{C}^\infty(M) : \mbox{$\int_M f=0$} \}.
\]
Precisely, the restriction $\int_M f=0$ means that the variation associated to $f$ is volume preserving.
The weak stability index of the hypersurface, denoted here by $\mathrm{Ind}_T(M)$, is characterized by
\[
\mathrm{Ind}_T(M)=
\max\{ \mathrm{dim}V : V\leqslant\mathcal{C}_T^\infty(M), \quad Q(f)<0 \quad \mbox{for every } f\in V \},
\]
and $M$ is called weakly stable if and only if $\mathrm{Ind}_T(M)=0$ (see \cite{A} for further details).

In \cite{BdCE}, Barbosa, do Carmo and Eschenburg characterized the totally umbilical spheres as the only compact
weakly stable constant mean curvature hypersurfaces in \s{n+1}. In \cite{ABP} the authors have recently showed
that the weak index of a compact constant mean curvature hypersurface $M^n$ in \s{n+1} which
is not totally umbilical and has constant scalar curvature is greater than or equal to $n+2$,
with equality if and only if $M^n$ is a Clifford hypersurface
$M_k(r)=\mathbb{S}^{k}(r)\times\mathbb{S}^{n-k}(\sqrt{1-r^2})$ with radius
$\sqrt{k/(n+2)}\leqslant r\leqslant\sqrt{(k+2)/(n+2)}$. Here we will complement this result by showing the following.
\begin{theorem}
\label{theorindex}
Let $M^n$ be a compact orientable hypersurface immersed into the Euclidean sphere
\s{n+1} with constant mean curvature. If $M$ has constant scalar curvature
and $M$ is neither a Clifford nor an umbilical hypersurface, then the weak
stability index of $M$ is greater than or equal to $2n+4$.
\end{theorem}
\begin{proof}
The condition on the scalar curvature implies that, $\|A\|^2$ is constant. Let us first consider the case where $H=0$.
Since $M^n$ is not totally umbilical (i.e., totally geodesic), then $\|A\|^2>0$. Even more, since $M$ is not a
minimal Clifford hypersurface we have that $\|A\|^2>n$, by a classical result due to \cite{S} and \cite{CDK,L} (see
\cite[Theorem 6]{A}). By Proposition\rl{prop5} we have that the functions $\ell_v$ and $f_v$ are eigenfunctions of the
Laplacian with positive eigenvalues $n$ and $\|A\|^2>n$, respectively (observe that with our criterion, a real
number $\lambda$ is an eigenvalue of $\Delta$ if and only if $\Delta u+\lambda u=0$ for some smooth function $u\in\mathcal{C}^\infty(M)$, $u\not\equiv 0$).
In particular, the functions $\ell_v$ and $f_v$ satisfy the condition $\int_M f=0$, and they also satisfy $J(\ell_v)=\|A\|^2\ell_v$ and
$Jf_v=nf_v$. That is, they are also eigenfunctions of $J$ with negative eigenvalues $-\|A\|^2$ and $-n$, respectively.
Let
\[
V_1=\{\ell_v : v\in\R{n+2}\} \quad\hbox{and} \quad V_2=\{f_v : v\in\R{n+2}\}.
\]
Then,
\beq
\label{indT1}
\mathrm{Ind}_T(M)\geq \mathrm{dim}(V_1\oplus V_2)=\mathrm{dim}V_1+\mathrm{dim}V_2,
\eeq
where the last equality is due to the fact that $V_1$ and $V_2$ are $L^2$-orthogonal subspaces, because they are
eigenspaces of $\Delta$ associated to different eigenvalues. Finally, as pointed out in Subsection\rl{ss1.2},
we also know that if either $\mathrm{dim}V_1<n+2$ or $\mathrm{dim}V_2<n+2$, then $M$ must be a totally geodesic sphere
(see \cite[Lemma 3.1]{P}). Therefore, in our case we have $\mathrm{dim}V_1=\mathrm{dim}V_2=n+2$, and by \rf{indT1} we
conclude that $\mathrm{Ind}_T(M)\geq 2n+4$.

We will now consider the case $H\ne 0$. By Cauchy-Schwarz inequality we
have that $\|A\|^2\geq nH^2$, and equality only occurs if $M$ is totally umbilical. In this case, following our ideas
in \cite{ABP}, we will work with test functions of the form $\ell_v-\alpha_{{\pm}}f_v$, where
$$
\alpha_{\pm}=\frac{\|A\|^2-n\pm\sqrt{D}}{2nH} \quad \hbox{with} \quad
D=(\|A\|^2-n)^2+4n^2H^2>0.
$$
Let
\[
U_{+}=\{\ell_v-\alpha_{+}f_v: v\in\R{n+2}\} \quad\hbox{and} \quad U_{-}=\{\ell_v-\alpha_{-}f_v: v\in\R{n+2}\}.
\]
Then, by Proposition\rl{prop5} we have that $\Delta u+\mu_{\pm}u=0$ for every $u\in U_{\pm}$, where
\[
0<\mu_{-}=\frac{n+\|A\|^2-\sqrt{D}}{2}<\mu_{+}=\frac{n+\|A\|^2+\sqrt{D}}{2},
\]
and, therefore, $Ju+\lambda_{\pm}u=0$ for every $u\in U_{\pm}$, with
\[
\lambda_{-}=\frac{-(n+\|A\|^2)-\sqrt{D}}{2}<\lambda_{+}=\frac{-(n+\|A\|^2)+\sqrt{D}}{2}<0
\]
(for the details, see \cite[Section 4]{ABP}). In particular, functions belonging to $U_{\pm}$ also satisfy the condition
$\int_M f=0$, and
\beq
\label{indT2}
\mathrm{Ind}_T(M)\geq \mathrm{dim}(U_{+}\oplus U_{-})=\mathrm{dim}U_{+}+\mathrm{dim}U_{-}.
\eeq
Finally, since $M$ is neither a totally umbilical sphere nor a Clifford hypersurface, our Theorem\rl{maintheorem}
implies that $\mathrm{dim}U_{+}=\mathrm{dim}U_{-}=n+2$, and by \rf{indT2} we
conclude that $\mathrm{Ind}_T(M)\geq 2n+4$. \eop

\end{proof}

\section*{Acknowledgements}
The authors would like to thank to the referee for valuable suggestions which improved the paper.

\end{document}